\documentclass[journal]{IEEEtran}

\usepackage{booktabs}
\usepackage{amsfonts}
\usepackage{mathrsfs} 
\usepackage{boldline}
\usepackage[table]{xcolor}
\usepackage{amsmath} 
\usepackage{psfrag} 
\usepackage{multirow} 
\usepackage{arydshln}
\usepackage{enumerate}
\usepackage{bbm}
\usepackage{graphicx}
\usepackage{color}
\usepackage{setspace}
\usepackage{mathtools}
\usepackage{enumitem}
\usepackage{algorithmic}
\usepackage{algorithm}
\usepackage{tabularx}
\usepackage{xfrac}
\usepackage{booktabs,multirow}
\newcommand{\revise}[1]{#1}

\newboolean{showcomments}
\setboolean{showcomments}{true}
\newcommand{\fz}[1]{\ifthenelse{\boolean{showcomments}}
{ \textcolor{blue}{(Fengyu:  #1)}}{}}
\newcommand{\slow}[1]{\ifthenelse{\boolean{showcomments}}
{ \textcolor{red}{(Steven:  #1)}}{}}
\newcommand{\ahmed}[1]{\ifthenelse{\boolean{showcomments}}
{ \textcolor{orange}{(Ahmed:  #1)}}{}}
\newcommand{\update}[1]{#1}

\newcommand{\Real}{\mathbb{R}}
\newcommand{\Complex}{\mathbb{C}}

\newcommand{\iu}{\mathrm{i}}

\newcommand{\T}{\mathsf{T}}
\newcommand{\Hn}{\mathsf{H}}

\newcommand{\rank}{{\rm rank}}
\newcommand{\diag}{{\rm diag}}

\newcommand{\tr}{{\rm tr}}

\newcommand{\zero}{\mathbf{0}}
\newcommand{\one}{\mathbf{1}}
\newcommand{\Mat}{\mathbf{M}}
\newcommand{\MWXrho}{\Mat^{\W,\X,\CurDeltaSqr}}
\newcommand{\MvSl}{\Mat^{\VolSqr,\PF,\CurSqr}}
\newcommand{\MvXrho}{\Mat^{\VolSqr,\X,\CurDeltaSqr}}
\newcommand{\AMat}{\mathbf{A}}
\newcommand{\BMat}{\mathbf{B}}
\newcommand{\CMat}{\mathbf{C}}

\newcommand{\wVec}{\mathbf{w}}
\newcommand{\xVec}{\mathbf{x}}
\newcommand{\yVec}{\mathbf{y}}
\newcommand{\zVec}{\mathbf{z}}

\newcommand{\Graph}{\mathcal{G}}
\newcommand{\Vertex}{\mathcal{V}}
\newcommand{\Edge}{\mathcal{E}}
\newcommand{\Nbus}{N}
\newcommand{\Y}{{\rm Y}}
\newcommand{\V}{\mathbf{V}}
\newcommand{\Vref}{\mathbf{V}_{\rm ref}}
\newcommand{\s}{\mathbf{s}}
\newcommand{\py}{\mathbf{p}}
\newcommand{\qy}{\mathbf{q}}
\newcommand{\pd}{\mathbf{p}_{\Delta}}
\newcommand{\qd}{\mathbf{q}_{\Delta}}
\newcommand{\pya}[1]{\mathbf{p}_{#1}}
\newcommand{\qya}[1]{\mathbf{q}_{#1}}
\newcommand{\pda}[1]{\mathbf{p}_{\Delta,#1}}
\newcommand{\qda}[1]{\mathbf{q}_{\Delta,#1}}
\newcommand{\pyref}{\mathbf{\overline{p}}_{j}^\phi}
\newcommand{\qyref}{\mathbf{\overline{q}}_{j}^\phi}
\newcommand{\pdref}{\mathbf{\overline{p}}_{\Delta,j}^\phi}
\newcommand{\qdref}{\mathbf{\overline{q}}_{\Delta,j}^\phi}
\newcommand{\phases}{\mathsf{\Phi}}

\newcommand{\Cur}{{\bf{I}}}
\newcommand{\X}{{\bf{X}}}
\newcommand{\sol}{{\bf{u}}}

\newcommand{\CurDeltaSqrAux}{\tilde{\rho}}
\newcommand{\MatDiff}{\Gamma}
\newcommand{\CurDeltaSqr}{\rho}
\newcommand{\W}{\mathbf{W}}
\newcommand{\PF}{\mathbf{S}}
\newcommand{\VolSqr}{\mathbf{v}}
\newcommand{\CurSqr}{\mathbf{\ell}}

\newcommand{\BIM}{{\rm BIM}}
\newcommand{\BFM}{{\rm BFM}}

\newcommand{\set}{\mathcal{S}}
\newcommand{\feasiblesetX}{\mathcal{X}}
\newcommand{\feasiblesetY}{\mathcal{Y}}

\newtheorem{assumption}{Assumption}
\newtheorem{proposition}{Proposition}
\newtheorem{definition}{Definition}
\newtheorem{remark}{Remark}

\newtheorem{lemma}{Lemma}
\newtheorem{theorem}{Theorem}

\begin{document}
%
\title{Exactness of OPF Relaxation on Three-phase Radial Networks with Delta Connections}
%
%
%

\author{Fengyu~Zhou,~\IEEEmembership{Student Member,~IEEE,}
        Ahmed~S.~Zamzam,~\IEEEmembership{Member,~IEEE,}
        Steven~H.~Low,~\IEEEmembership{Fellow,~IEEE,}
        and~Nicholas~D.~Sidiropoulos,~\IEEEmembership{Fellow,~IEEE}
\thanks{
{
}

Fengyu Zhou and Steven H. Low are with the Department of Electrical Engineering, California Institute of Technology, Pasadena,
CA 91125 USA (e-mail: \{f.zhou, slow\}@caltech.edu). Ahmed S. Zamzam is with the National Renewable Energy Laboratory, Golden, CO 80401 USA (e-mail: ahmed.zamzam@nrel.gov). Nicholas D. Sidiropoulos is with the Electrical and Computer Engineering Department, University of Virginia, Charlottesville, VA 22901 USA (email: nikos@virginia.edu)\newline This work was authored in part by the National Renewable Energy Laboratory, operated by Alliance for Sustainable Energy, LLC, for the U.S. Department of Energy (DOE) under Contract No. DE-AC36-08GO28308. The work of A. S. Zamzam was supported by the
Laboratory Directed Research and Development Program at the National
Renewable Energy Laboratory.  The views expressed in the article do not necessarily represent the views of the DOE or the U.S. Government. The U.S. Government retains and the publisher, by accepting the article for publication, acknowledges that the U.S. Government retains a nonexclusive, paid-up, irrevocable, worldwide license to publish or reproduce the published form of this work, or allow others to do so, for U.S. Government purposes.}
}

\markboth{}
{Shell \MakeLowercase{\textit{et al.}}: Bare Demo of IEEEtran.cls for IEEE Journals}
%



\maketitle

\begin{abstract}
Simulations have shown that while semi-definite relaxations of AC optimal power flow (AC-OPF) on three-phase radial networks with only wye connections tend to be exact, the presence of delta connections seem to render them inexact.
\revise{This paper shows that the such inexactness originates from the non-uniqueness of relaxation solutions and numerical errors amplified by the non-uniqueness.
This finding motivates two algorithms to recover the exact solution of AC-OPF in the presence of delta connections.} 
In simulations using IEEE 13, 37 \update{and 123-bus} systems, the proposed algorithms provide exact optimal solutions up to numerical precision.
\end{abstract}


%
\IEEEpeerreviewmaketitle

\section*{Nomenclature}
\addcontentsline{toc}{section}{Nomenclature}
\begin{IEEEdescription}[\IEEEusemathlabelsep\IEEEsetlabelwidth{$V_1,V_2,V_3$}]
\update{
\item[$\Graph(\Vertex,\Edge)$] Graph with vertices $\Vertex$ and edges $\Edge$.
\item[$z_{jk},y_{jk}$] Impedance and admittance matrices for line $(j,k)$.
\item[$\V_j,\V_j^\phi$] Complex voltage at bus $j$ (phase $\phi$).
\item[$\underline{\V},\overline{\V}$] Limits on voltage magnitude, which may be set as $\underline{\V}=\underline{V}\one$ and $\overline{\V}=\overline{V}\one$.
\item[$\s_j, \s_j^\phi$] Wye injection at bus $j$ (phase $\phi$).
\item[$\py_j, \qy_j$] Active and reactive wye injection at bus $j$, i.e., $\s_j=\py_j+\iu\qy_j$.
\item[$\s_{\Delta,j}, \s_{\Delta,j}^\phi$] Delta injection at bus $j$ (phase $\phi$).
\item[$\py_{\Delta,j}, \qy_{\Delta,j}$] Active and reactive delta injection at bus $j$, i.e., $\s_{\Delta,j}=\py_{\Delta,j}+\iu\qy_{\Delta,j}$.
\item[$f(\s,\s_{\Delta})$] Cost function of AC-OPF.
\item[$\set_j, \set_{\Delta,j}$] Feasible sets for $\s_j$ and $\s_{\Delta,j}$.
\item[$\Cur_{jk}$] Sending-end currents from $j$ to $k$.
\item[$\Cur_{\Delta,j}$] Delta line currents at bus $j$.
\item[$\PF_{jk}$] Sending-end branch power flow, i.e., $\V_{j}\Cur_{jk}^\Hn$.
\item[$\X_j$] Power flow on delta devices, i.e., $\V_{j}\Cur_{\Delta,j}^\Hn$.
\item[$\W$,$\VolSqr_j$] $\W=\V\V^\Hn$, $\VolSqr_j=\V_j\V_j^\Hn$.
\item[$\CurSqr_{jk}, \CurDeltaSqr_j$]$\CurSqr_{jk}=\Cur_{jk}\Cur_{jk}^\Hn$, $\CurDeltaSqr_j=\Cur_{\Delta,j}\Cur_{\Delta,j}^\Hn$.
\item[$\MWXrho$] Block matrix defined in \eqref{eq:WXrho}.
\item[$\MvSl$] Block matrix defined in \eqref{eq:vSl}.
\item[$\MvXrho$] Block matrix defined in \eqref{eq:vXrho}.
\item[$p_{\rm loss}$] Total active power loss in the network.
\item[$d_p$, $d_q$] Penalty functions on active and reactive power deviation.
\item[${\phases}_{\Y}^j,{\phases}_{\Delta}^j$] Available wye and delta connected phases at bus $j$.}
\end{IEEEdescription}

\section{Introduction}\label{sec:intro}
\IEEEPARstart{O}{ptimal}
 power flow (OPF) is a mathematical program that finds an optimal operating point for a power grid subject to laws of physics and operational constraints~\cite{carpentier1962contribution}.
 OPF formulated under the AC model is known to be both nonconvex and NP-hard to solve~\cite{bienstock2019strong, lehmann2016ac}.
Commonly used methods to approximately solve this non-convex problem include, e.g., the Newton-Raphson method~\cite{momoh1999review} and various linearizations~\cite{wells1968method,burchett1982large,Bernstein2018load}. 
Another approach that has emerged over the last decade or so is to relax the OPF problem to a convex program, such as semi-definite program (SDP) relaxations and second-order cone program (SOCP) relaxations \cite{jabr2006radial, bai2008semidefinite, farivar2013branch}. 
For single-phase radial networks (i.e., networks with a tree topology) as well as the single-phase equivalent of a balanced three-phase radial network, simulations have shown that these relaxations often yield solutions that are also global optima of the original nonconvex problems \cite{farivar2012,kekatos2014}. Sufficient conditions that guarantee exact relaxations for single-phase radial networks have subsequently been proved.
There is now sizeable literature on OPF relaxations; see, e.g., comprehensive surveys in \cite{molzahn2019survey, low2014convexII} for pointers to various contributions including many earlier surveys on OPF.
Many recent works further refine and improve the convexification of OPF problems and exhibit promising performance \cite{zamzam2018beyond, zamzam2017qcqp}. 

Most of this literature focuses on single-phase models, but distribution systems have multiple phases that are increasingly unbalanced as distributed energy resources continue to grow.
Mathematically, we can identify each bus-phase pair of a three-phase radial network with an equivalent single-phase bus and transform it into a single-phase equivalent circuit with a meshed topology \cite{berg1967mechanized,laughton1968analysis}.
\update{Most analytical results on exact relaxation for meshed networks are restricted to weakly-cyclic networks \cite{madani2014convex},
but the single-phase equivalent of a three-phase radial network is beyond this class even for a 2-bus toy example.}
SDP relaxation has been generalized in~\cite{dall2013distributed,gan2014convex} to three-phase radial networks with only wye connections and shown to be exact in the simulation of several test cases.
A recent work \cite{zhou2019sufficient} proves a sufficient condition for exact relaxation in this case.  
\update{
Besides SDP relaxation, \cite{wang2017chordal} also provides an iterative algorithm for three-phase networks, also without delta connections.}


\update{
Semi-definite relaxation is recently extended in \cite{zhao2017optimal} to networks with delta connected devices by introducing a new positive semi-definite matrix that represents the outer product of voltages and phase-to-phase currents in the delta connections (matrix $\MvXrho(j)$ in \eqref{eq:vXrho} below). 
Simulation results in~\cite{zhao2017optimal} 
showed that, surprisingly, 
this matrix was never rank-1 at an optimal solution of the relaxation. This seems to suggest that the SDP relaxation was inexact in these simulations.
In this paper, we show that even though the matrix $\MvXrho(j)$ fails to attain rank 1, an exact solution can still be recovered under certain conditions; see Theorem \ref{thm:main-BFM} and Remark \ref{rmk:numerical}. The inexactness in previous works is due to two issues.
First, optimal solutions to the SDP relaxation in these simulations are generally not unique, and the exact solution is only one of them which is not returned by the solver.
Second, such non-uniqueness could significantly amplify the numerical error and make it computationally challenging to recover the exact solution.
We propose two variants of the standard SDP relaxation that address both issues.
The first algorithm applies post-processing to the relaxation solution and tends to provide lower cost but larger constraint violation, while the second algorithm adds a penalty term to the cost and tends to provide higher cost but smaller constraint violation.
Simulations of both algorithms corroborate the theoretical results and show that they can recover exact solutions for three IEEE distribution feeders.
}


The remainder of the paper is organized as follows. In Section \ref{sec:model}, we define the network structure and formulate the three-phase OPF problem in both the bus injection model (BIM) and the branch flow model (BFM). We subsequently show that the two models are essentially equivalent to each other. In Section \ref{Sec:algorithm}, we prove that the global optimal solution to the nonconvex OPF problem can be recovered from its relaxation under certain conditions, and two algorithms are presented. Section \ref{sec:analysis} shows the equivalence between BIM and BFM.
Finally, in Section \ref{sec:simulation}, we apply our algorithms to IEEE 13-, 37- \update{and 123-bus} systems.
\section{System Model}\label{sec:model}
\subsection{Network Structure}
\begin{figure}
\centering
\update{
\includegraphics[width=\columnwidth]{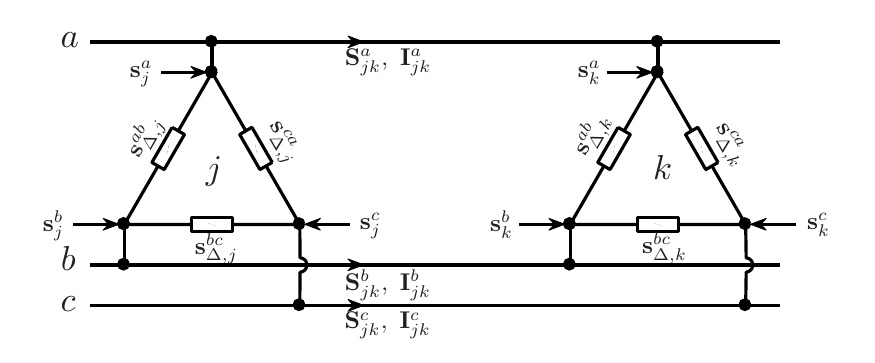}
\caption{Illustration of line $(j,k)$ with both wye and delta connections.}
\label{fig:line}}
\end{figure}

We study the model proposed in ~\cite{gan2014convex, zhao2017optimal}.
Let the directed graph representing the electrical network be $\Graph=(\Vertex,\Edge)$, where $\Vertex=\{0,1,\dots,n\}$ denotes the set of buses, and $\Edge\subseteq\Vertex\times\Vertex$ denotes the set of edges; and let $\Nbus:=|\Vertex|=n+1$. 
In this paper, we focus on the case where $\Graph$ represents a radial network (i.e., a tree) because most distribution networks have a tree topology.
Throughout the paper, we will use (graph, vertex, edge) and (power network, bus, line) interchangeably. 
Without loss of generality, we let bus $0$ be the substation bus where the distribution feeder is connected to a transmission network. 
Suppose the substation also serves as the slack bus, so the voltages at the substation bus are fixed and specified.
We use $j\rightarrow k$ to denote a directed edge from bus $j$ to $k$.
In many situations, when we do not care about the direction of the edge,
we simply use $(j,k)$ and $j\sim k$ interchangeably to denote an edge connecting bus $j$ and $k$.
That means either $j\rightarrow k$ or $k\rightarrow j$ is in $\Edge$.
Consider a three-phase line $(j,k)$ characterized by the series impedance matrix
$z_{jk} \in \Complex^{3\times 3}$.
When line $(j, k)$ has three phases, the inverse of $z_{jk}$, denoted as $y_{jk}$, is the admittance of line $(j, k)$. If branch $(j,k)$ has less than three phases, then we fill the rows and columns of $z_{jk}$ corresponding to the missing phases with zeros, and we let the admittance matrix $y_{jk}$ be the pseudo-inverse of $z_{jk}$.
Last, let $y_j\in\Complex^{3\times 3}$ denote the admittance of a shunt device connected to bus $j$.\footnote{The shunt here refers to a capacitive device at bus $j$ and not the line charging in the $\Pi$ circuit model.}

For each bus $j$, let the voltages of all three phases at bus $j$ be collected in the vector $\V_j\in\Complex^3$. We use $\V_j^{\phi}$ for $\phi\in\{a,b,c\}$ to indicate the voltage of phase $\phi$.
The voltage $\V_0$ at slack bus $0$ is known and denoted by $\Vref$.
Let $\V=[\V_0^\T,\V_1^\T,\dots,\V_{n}^\T]^\T$ collect the voltages for the entire network. Similarly, we use $\s_j^\phi$ to denote the bus injection for phase $\phi$ at bus $j$; and we denote $\s_j$ and $\s$ as the injections at bus $j$ and in the entire network, respectively.

For delta connected components, we use $\Cur_{\Delta,j}\in\Complex^3$ to collect the delta line currents for phases in $\{ab,bc,ca\}$. 
Define the matrix
\begin{align*}
\MatDiff:=\left[
\begin{array}{ccc}
1 & -1 & 0\\
0 & 1 & -1\\
-1 & 0 & 1
\end{array}
\right].
\end{align*}
Therefore, the complex power injections of the delta connected components at bus $j$ can be expressed as $\s_{\Delta,j}=\diag(\MatDiff \V_j\Cur_{\Delta,j}^\Hn)$. 
The net nodal injections
contributed by delta connections at bus $j$ are given by $-\diag(\V_j\Cur_{\Delta,j}^\Hn\MatDiff)$  \update{(see the illustration in Fig. \ref{fig:line}).}
Assume that the operation regions for $\s_j$ and $\s_{\Delta,j}$ at bus $j$ are convex compact sets $\set_j$ and $\set_{\Delta,j}$, respectively.

The AC power flow equations are:
\begin{subequations}
\begin{align}
\nonumber
&\s_j-\diag(\V_j\Cur_{\Delta,j}^\Hn\MatDiff)-\diag(\V_j\V_j^\Hn y_j^\Hn)\\
=&\sum\limits_{k:j\sim k}\diag((\V_j\V_j^\Hn-\V_j\V_k^\Hn)y_{jk}^\Hn)
\label{eq:pf.a}\\
&\s_{\Delta,j} = \diag(\MatDiff\V_j\Cur_{\Delta,j}^\Hn).
\label{eq:pf.b}
\end{align}\label{eq:pf}
\end{subequations}
where \eqref{eq:pf.a} is the power balance equation at bus $j$ and \eqref{eq:pf.b} defines the power through delta connected components. 

Similar to \cite{zhao2017optimal}, we will adopt $\X_{j},\CurDeltaSqr_j\in\Complex^{3\times 3}$ as auxiliary matrices in the following models as well as their corresponding relaxations. The outer products of voltages and currents are defined as follows:
\begin{subequations}
\begin{align}
\X_{j}=\V_j\Cur_{\Delta,j}^\Hn
\label{eq:X}\\
\CurDeltaSqr_j=\Cur_{\Delta,j}\Cur_{\Delta,j}^\Hn.
\label{eq:rho}
\end{align}
\end{subequations}

We consider an OPF problem that minimizes a cost function $f(\s,\s_{\Delta})$ over variables $(\s,\s_{\Delta},\V,\Cur_{\Delta})$ subject to power flow equations \eqref{eq:pf} as well as voltage and injection limits:
\begin{subequations}
\begin{align}
\underset{\s,\s_\Delta,\V,\Cur_{\Delta}}{\text{minimize}}  &&& f(\s,\s_{\Delta})
\label{eq:OPF.a}\\
\text{\quad subject to}&&& \eqref{eq:pf}
\label{eq:OPF.b}\\
&&& \V_0 = \Vref
\label{eq:OPF.c}\\
&&&\s_j\in\set_j,~\s_{\Delta,j}\in\set_{\Delta,j},~\text{for $j\in\Vertex$}
\label{eq:OPF.d}\\
&&&\underline{\V}\leq|\V|\leq\overline{\V}.
\label{eq:OPF.e}
\end{align}
\label{eq:OPF}
\end{subequations}
In \eqref{eq:OPF.e}, $|\V|$ stands for the modulus of $\V$ elementwise, and $\underline{\V},\overline{\V}\in\Real^{3N}$ are the lower and upper limits for the voltages. If the limits are homogeneous across all buses and phases, we can denote them as $\underline{V}\one, \overline{V}\one$, where $\underline{V},\overline{V}$ are scalars and $\one$ is the all-one vector.

\subsection{Bus Injection Model}
The {\it bus injection model} (BIM) is defined in terms of $(\s,\s_\Delta,\W,\X,\CurDeltaSqr)$, where we use $\W\in\Complex^{3\Nbus\times 3\Nbus}$ to replace $\V\V^\Hn$ in \eqref{eq:pf}.
The matrix $\W_{jk}\in\Complex^{3\times 3}$ is the $(j,k)$ submatrix of $\W$. 
For notational simplicity, we let
\begin{align}\label{eq:WXrho}
\MWXrho(j):=&\left[
\begin{array}{cc}
\W_{jj} & \X_j \\
\X_j^\Hn & \CurDeltaSqr_j
\end{array}
\right] &&\text{for $j\in\Vertex$}.
\end{align}
The power flow model is represented as
\begin{subequations}
\begin{align}
\nonumber
&\s_j-\diag(\X_j\MatDiff)-\diag(\W_{jj} y_j^\Hn)\\
=&\sum\limits_{k:j\sim k}\diag((\W_{jj}-\W_{jk})y_{jk}^\Hn)
\label{eq:PF-BIM.a}\\
&\s_{\Delta,j} = \diag(\MatDiff\X_j)
\label{eq:PF-BIM.b}\\
&\W_{00}=\Vref\Vref^\Hn
\label{eq:PF-BIM.c}\\
&\W \succeq 0
\label{eq:PF-BIM.d}\\
&\rank(\W)=1
\label{eq:PF-BIM.e}\\
&\MWXrho(j) \succeq 0
\label{eq:PF-BIM.f}\\
&\rank\left(\MWXrho(j)\right)=1.
\label{eq:PF-BIM.g}
\end{align}
\label{eq:PF-BIM}
\end{subequations}

Hence, the AC-OPF problem in BIM formulation is: 
\begin{subequations}
\begin{align}
\underset{\s,\s_\Delta,\W,\X,\CurDeltaSqr}{\text{minimize}}  &&& f(\s,\s_{\Delta})
\label{eq:OPF-BIM.a}\\
\text{\quad subject to}&&& \eqref{eq:PF-BIM}, \eqref{eq:OPF.d}.
\label{eq:OPF-BIM.b}\\
&&& \diag(\underline{\V}\underline{\V}^\Hn)\leq\diag(\W)\leq\diag(\overline{\V}\overline{\V}^\Hn)
\label{eq:OPF-BIM.c}
\end{align}
\label{eq:OPF-BIM}
\end{subequations}

\subsection{Branch Flow Model}
In a {\it branch flow model} (BFM), we introduce $\PF$, $\VolSqr$, and $\CurSqr$ to model the branch power flow, squared voltages, and squared currents, respectively. The matrices $\PF$, $\VolSqr$, and $\CurSqr$ can be written as
\begin{subequations}
\begin{align}
&\PF=(\PF_{jk}\in\Complex^{3\times 3})_{(j\rightarrow k)\in\Edge},&&\PF_{jk}=\V_j\Cur_{jk}^\Hn\\
&\VolSqr=(\VolSqr_j\in\Complex^{3\times 3})_{j\in\Vertex},&&\VolSqr_j=\V_j\V_j^\Hn\\
&\CurSqr=(\CurSqr_{jk}\in\Complex^{3\times 3})_{(j\rightarrow k)\in\Edge},&&\CurSqr_{jk}=\Cur_{jk}\Cur_{jk}^\Hn.
\end{align}
\end{subequations}
Here, $\Cur_{jk}:=y_{jk}(\V_j-\V_k)$ is the sending-end current from bus $j$ to bus $k$, and $\PF_{jk}$ is the sending-end branch power from $j$ to $k$.
Let
\begin{subequations}
\begin{align}
\MvSl(j,k):=&\left[
\begin{array}{cc}
\VolSqr_j & \PF_{jk} \\
\PF_{jk}^\Hn & \CurSqr_{jk}
\end{array}
\right] &&\text{for $j\rightarrow k$}\label{eq:vSl}\\
\MvXrho(j):=&\left[
\begin{array}{cc}
\VolSqr_j & \X_j \\
\X_j^\Hn & \CurDeltaSqr_j
\end{array}
\right] &&\text{for $j\in\Vertex$}.\label{eq:vXrho}
\end{align}
\end{subequations}
The branch flow model is defined in terms of $(\s,\s_\Delta,\PF,\VolSqr,\CurSqr,\X,\CurDeltaSqr)$, and it is expressed as
\begin{subequations}
\begin{align}
&\VolSqr_k=\VolSqr_j-(\PF_{jk}z_{jk}^\Hn+z_{jk}\PF_{jk}^\Hn)+z_{jk}\CurSqr_{jk}z_{jk}^\Hn
\label{eq:PF-BFM.a}\\
\begin{split}
&\sum\limits_{k:j\rightarrow k}\diag(\PF_{jk})
-\sum\limits_{l:l\rightarrow j}\diag(\PF_{lj}-z_{lj}\CurSqr_{lj})\\
\qquad\qquad=&-\diag(\VolSqr_j y_j^\Hn+\X_j\MatDiff)+\s_j
\end{split}
\label{eq:PF-BFM.b}\\
&\s_{\Delta,j}=\diag(\MatDiff\X_j)
\label{eq:PF-BFM.c}\\
& \VolSqr_0 = \Vref\Vref^\Hn
\label{eq:PF-BFM.d}\\
&\MvSl(j,k) \succeq 0
\label{eq:PF-BFM.g}\\
&\rank\left(\MvSl(j,k)\right)= 1
\label{eq:PF-BFM.h}\\
&\MvXrho(j) \succeq 0
\label{eq:PF-BFM.i}\\
&\rank\left(\MvXrho(j)\right)= 1.
\label{eq:PF-BFM.j}
\end{align}
\label{eq:PF-BFM}
\end{subequations}

The AC-OPF problem in the BFM form can be formulated as:
\begin{subequations}
\begin{align}
\underset{\s,\s_\Delta,\PF,\VolSqr,\CurSqr,\X,\CurDeltaSqr}{\text{minimize}} &&& f(\s,\s_{\Delta})
\label{eq:OPF-BFM.a}\\
\text{\quad subject to}&&& \eqref{eq:PF-BFM}, \eqref{eq:OPF.d}
\label{eq:OPF-BFM.b}\\
&&&\!\! \diag(\underline{\V}_j\underline{\V}_j^\Hn)\!\leq\!\diag(\VolSqr_j)\leq\diag(\overline{\V}_{\!j}\overline{\V}_{\!j}^\Hn)
\label{eq:OPF-BFM.c}
\end{align}
\label{eq:OPF-BFM}
\end{subequations}

\section{Main Results}\label{Sec:algorithm}
The main challenge to solving OPF problems \eqref{eq:OPF-BIM} and \eqref{eq:OPF-BFM} is the nonconvex rank constraints in \eqref{eq:PF-BIM.e}, \eqref{eq:PF-BIM.g}, \eqref{eq:PF-BFM.h} and \eqref{eq:PF-BFM.j}. 
If we drop all the rank-1 constraints, then we obtain
\begin{subequations}
\begin{align}
\underset{\s,\s_\Delta,\W,\X,\CurDeltaSqr}{\text{minimize}}  &&& f(\s,\s_{\Delta})
\label{eq:ROPF-BIM.a}\\
\text{\quad subject to}&&& \eqref{eq:PF-BIM.a}-\eqref{eq:PF-BIM.d}, \eqref{eq:PF-BIM.f},
\eqref{eq:OPF.d}, \eqref{eq:OPF-BIM.c}.
\label{eq:ROPF-BIM.b}
\end{align}
\label{eq:ROPF-BIM}
\end{subequations}
as the relaxation for the BIM and
\begin{subequations}
\begin{align}
\underset{\s,\s_\Delta,\PF,\VolSqr,\CurSqr,\X,\CurDeltaSqr}{\text{minimize}} &&& f(\s,\s_{\Delta})
\label{eq:ROPF-BFM.a}\\
\text{\quad subject to}&&& \eqref{eq:OPF.d}, \eqref{eq:OPF-BFM.c}, \eqref{eq:PF-BFM.a}-\eqref{eq:PF-BFM.g}, \eqref{eq:PF-BFM.i}
\label{eq:ROPF-BFM.b}
\end{align}
\label{eq:ROPF-BFM}
\end{subequations}
as the relaxation for the BFM.
Solving the relaxed problems~\eqref{eq:ROPF-BIM} and~\eqref{eq:ROPF-BFM} could lead to solutions that are infeasible for the original nonconvex problems~\eqref{eq:OPF-BIM} and~\eqref{eq:OPF-BFM} respectively when the solutions do not satisfy the rank-1 constraints.
In what follows, we will explore conditions under which optimal solutions of~\eqref{eq:OPF-BIM} and~\eqref{eq:OPF-BFM} can be recovered from their respective relaxations. First, the following lemma is presented, which is the main ingredient in the proofs of subsequent results.

\begin{lemma}\label{Lem:rank1}
Consider a block Hermitian matrix
\begin{align}\label{eq:block-psd}
\Mat:=
\left[
\begin{array}{cc}
    \AMat & \BMat \\
    \BMat^\Hn & \CMat
\end{array}
\right]
\end{align}
where $\AMat$ and $\CMat$ are both square matrices. 
If $\Mat\succeq 0$ and $\AMat=\xVec\xVec^\Hn$ for some vector $\xVec$, 
then there must exist some vector $\yVec$ such that $\BMat = \xVec\yVec^\Hn$. 
\end{lemma}
\begin{IEEEproof}
As $\Mat\succeq 0$, it can be decomposed as
\begin{align}\label{eq:psd-factorization}
\Mat=\left[
\begin{array}{cc}
    \Mat_1 \\
    \Mat_2
\end{array}
\right]
\left[
\begin{array}{cc}
    \Mat_1^\Hn & \Mat_2^\Hn
\end{array}
\right]
\end{align}
and $\AMat=\Mat_1\Mat_1^\Hn$, $\BMat=\Mat_1\Mat_2^\Hn$, $\CMat=\Mat_2\Mat_2^\Hn$.
Because $\AMat=\xVec\xVec^\Hn$ has rank-1, matrix $\Mat_1$ is in the column space of $\xVec$ and has rank-1 as well. There must exist vector $\zVec$ such that $\Mat_1=\xVec\zVec^\Hn$.
As a result, $\BMat=\Mat_1\Mat_2^\Hn=\xVec\zVec^\Hn\Mat_2^\Hn=\xVec(\Mat_2\zVec)^\Hn$.
\end{IEEEproof}

\update{
One observation in Lemma \ref{Lem:rank1} is when submatrices $\AMat$ and $\BMat$ are fixed and specified as $\xVec\xVec^\Hn$ and $\xVec\yVec^\Hn$, there are non-unique $\CMat$ to make $\Mat$ positive semi-definite.
Similarly, in the relaxations \eqref{eq:ROPF-BIM} and \eqref{eq:ROPF-BFM}, the optimal solutions are always non-unique. Taking \eqref{eq:ROPF-BIM} as an example, for any optimal solution $(\s^*,\s_{\Delta}^*,\W^*,\X^*,\CurDeltaSqr^*)$, one could add to $\CurDeltaSqr^*$ an arbitrary PSD matrix to obtain a different optimal solution $(\s^*,\s_{\Delta}^*,\W^*,\X^*,\CurDeltaSqr^*+\mathbf{K}\mathbf{K}^\Hn)$. This non-uniqueness in the optimal $\rho$ explains why in existing literature such as \cite{zhao2017optimal}, the relaxation \eqref{eq:ROPF-BIM} could compute rank-1 $\W$ (within numerical tolerance) but the resulting $\MWXrho$ is always not rank-1.
In fact, the next result shows in theory, if the optimal $\W$ is perfectly of rank 1 without any numerical error, then an feasible and optimal solution of \eqref{eq:OPF-BIM} is recoverable.}

\begin{theorem}\label{thm:main-BIM}
If $\sol^*=(\s^*,\s_{\Delta}^*,\W^*,\X^*,\CurDeltaSqr^*)$ is an optimal solution to \eqref{eq:ROPF-BIM} that satisfies $\rank(\W^*)=1$, then an optimal solution of \eqref{eq:OPF-BIM} can be recovered from $\sol^*$.
\end{theorem}
\begin{IEEEproof}
We decompose $\W_{jj}^*$ as $\V_j\V_j^\Hn$ for each $j$, where $\V_j$ is a vector. By Lemma \ref{Lem:rank1}, there exists vector $\Cur_{\Delta,j}$ such that $\X_j^*=\V_j\Cur_{\Delta,j}^\Hn$. 
One could construct $\tilde{\CurDeltaSqr}$ such that $\tilde{\CurDeltaSqr}_j=\Cur_{\Delta,j}\Cur_{\Delta,j}^\Hn$.

Since \eqref{eq:ROPF-BIM} is a relaxation of \eqref{eq:OPF-BIM},
for $(\s^*,\s_\Delta^*,\W^*,\X^*,\tilde{\CurDeltaSqr})$ to be optimal for \eqref{eq:OPF-BIM}, 
it is sufficient that it is feasible for \eqref{eq:OPF-BIM}.
Clearly, constraints \eqref{eq:OPF.d}, \eqref{eq:OPF-BIM.c}, \eqref{eq:PF-BIM.a}--\eqref{eq:PF-BIM.d} are satisfied because they are also the constraints in \eqref{eq:ROPF-BIM} and they do not involve the decision variable $\CurDeltaSqr$. Constraint \eqref{eq:PF-BIM.e} also holds as $\text{rank}(\W^*) = 1$.
Further, by Lemma \ref{Lem:rank1}, we have
\begin{align*}
\left[
\begin{array}{cc}
\W_{jj}^* & \X_j^* \\
(\X_j^*)^\Hn & \tilde{\CurDeltaSqr}_j
\end{array}
\right]=\left[
\begin{array}{c}
\V_j \\
\Cur_{\Delta,j}
\end{array}
\right]
\left[
\begin{array}{c}
\V_j \\
\Cur_{\Delta,j}
\end{array}
\right]^\Hn
\end{align*}
is both positive semi-definite and of rank-1. 
Hence, \eqref{eq:PF-BIM.f} and \eqref{eq:PF-BIM.g} are also satisfied.
Hence, $(\s^*,\s_\Delta^*,\W^*,\X^*,\tilde{\CurDeltaSqr})$ is feasible for \eqref{eq:OPF-BIM}, and this completes the proof.
\end{IEEEproof}

\begin{theorem}\label{thm:main-BFM}
If $\sol^*=(\s^*,\s_{\Delta}^*,\PF^*,\VolSqr^*,\CurSqr^*,\X^*,\CurDeltaSqr^*)$ is an optimal solution to \eqref{eq:ROPF-BFM} and satisfies $\rank(\Mat^{\VolSqr^*,\PF,^*\CurSqr^*}(j,k))=1$ for $j\sim k$ and $\rank(\VolSqr^*_j)=1$ for  $j\in\Vertex$, then an optimal solution of \eqref{eq:OPF-BFM} can be recovered from $\sol^*$.
\end{theorem}

The proof of Theorem \ref{thm:main-BFM} is omitted because it is similar to the proof of Theorem \ref{thm:main-BIM}.

\update{
Theorem \ref{thm:main-BIM} asserts that in theory, the only critical non-convex constraint of \eqref{eq:OPF-BIM} is \eqref{eq:PF-BIM.e}, in the sense that a solution satisfying \eqref{eq:PF-BIM.g} could always be recovered whenever \eqref{eq:PF-BIM.e} holds.
However in practice, $\W^*$ is typically not exactly rank-1 due to numerical precision and therefore Theorem \ref{thm:main-BFM} could not be directly applied to recover the optimal solution as long as numerical error exists.
This is because the recovery method in Theorem \ref{thm:main-BIM} relies on the rank-1 decomposition of $\X^*$. In practice even if $\W^*$ is close to being rank-1, the optimal $\X^*$ could still be very different from being rank-1, as we will explain below in Remark \ref{rmk:numerical}.

\begin{remark}[Spectrum Error] \label{rmk:numerical}
More specifically, the matrix $\AMat$ in~\eqref{eq:psd-factorization} being \emph{approximately} rank-1 does not necessarily mean that $\BMat$ is also approximately rank-1.\footnote{Here, being approximately rank-1 means that the second largest eigenvalue of the matrix is nonzero but smaller than the largest eigenvalue by several orders of magnitude.}
For example, consider the case where ${\bf x}$, ${\bf e}_1$, and ${\bf e}_2$ are orthogonal vectors with norms $1$, $10^{-4}$, and $10^{-5}$, respectively. Similarly, let ${\bf y}$, ${\bf z}_1$, and ${\bf z}_2$ be orthogonal vectors with norms $1$, $10^{4}$, and $10^5$, respectively. Then, construct the matrix ${\bf M}$ as in~\eqref{eq:psd-factorization} with $\Mat_1 = [{\bf x}\ \ {\bf e}_1\ \ {\bf e}_2]$ and $\Mat_2 = [{\bf y}\ \ {\bf z}_1\ \ {\bf z}_2]$. 
Clearly, the resulting matrix $\Mat$ has the upper left diagonal block that is approximately rank-1. In fact, the ratio between the two leading eigenvalues of the upper diagonal block is $10^{-8}$. On the other hand, the upper right block is of rank 3 with three singular values of $1$.
Consequently, even when $\W^*$ is close to rank-1 within a certain numerical tolerance, $\X^*$ could be far from being a rank-1 matrix, especially if $\CurDeltaSqr^*$ already contains a large redundant PSD matrix $\mathbf{K}\mathbf{K}^\Hn$.
Decomposing $\X^*$ as the product of two vectors, as in the proof of Theorem \ref{thm:main-BIM}, could result in a large numerical error.
In other words, the non-uniqueness of $\CurDeltaSqr$ could significantly amplify the numerical error.
In fact, as long as the second largest eigenvalue of $\W^*$ is not exactly $0$, then no matter how small it is, such spectrum error could potentially be significant, especially when the trace of $\mathbf{K}\mathbf{K}^\Hn$ is large.
\end{remark}

To summarize, there are two factors that prevent the relaxation output from being exact. The first is the non-uniqueness in the relaxation solution, and the second is that such non-uniqueness further greatly amplify the numerical error in computation.
This finding motivates two algorithms for practical implementation.}

\subsection{Relaxation with Post-Processing}
Remark \ref{rmk:numerical} shows even if $\W^*$ is approximately rank-1, $\X^*$ could be far away from being rank-1, and thus recovering the vector $\Cur_{\Delta,j}$ from $\X_j^*$ can lead to poor numerical performance. 
In the first algorithm, we instead recover $\Cur_{\Delta,j}^\Hn$ as $\diag^{-1} (\MatDiff \V_j) \s_{\Delta,j}^*$ from \eqref{eq:pf.b}, 
and then we reconstruct $\tilde{\X_j}$ as $\V_j\Cur_{\Delta,j}^\Hn$.
If there is no numerical error, $\X^*$ and $\tilde{\X}$ should be equal; however, in the presence of spectrum error, they could be different, as discussed in Remark \ref{rmk:numerical}.
The pseudo code is provided in Algorithm \ref{Alg:post-proc}.
\begin{algorithm}
\caption{Relaxation Algorithm with Post-Processing.}
\begin{algorithmic}[1]\label{Alg:post-proc}
 \renewcommand{\algorithmicrequire}{\textbf{Input:}}
 \renewcommand{\algorithmicensure}{\textbf{Output:}}
 \REQUIRE $y$, $\set$, $\set_{\Delta}$
 \ENSURE  Optimal solution $(\s,\s_\Delta,\W,\X,\CurDeltaSqr)$ to \eqref{eq:OPF-BIM}.
  \STATE Solve \eqref{eq:ROPF-BIM} to obtain $(\s^*,\s_\Delta^*,\W^*, \X^*, \CurDeltaSqr^*$).
  \IF {($\rank(\W^*)>1$)}
  \STATE Output `Failed!'
  \STATE Exit
  \ELSE
  \STATE Decompose $\W_{jj}^*=\V_j \V_j^\Hn$
  \STATE $ \Cur_{\Delta,j}^\Hn \leftarrow \text{diag}^{-1} (\MatDiff \V_j) \s_{\Delta,j}^*   $
  \STATE $\tilde{\X} \leftarrow \V_j \Cur_{\Delta,j}^\Hn$,
  $\CurDeltaSqrAux_j\leftarrow\Cur_{\Delta,j}\Cur_{\Delta,j}^\Hn$
  \RETURN $(\s^*,\s_\Delta^*,\W^*,\tilde{\X},\CurDeltaSqrAux)$ 
  \ENDIF
 \end{algorithmic} 
 \end{algorithm}


\begin{theorem}\label{thm:alg-post-proc}
If Algorithm \ref{Alg:post-proc} does not fail, then its output is an optimal solution of \eqref{eq:OPF-BIM}.
\end{theorem}

Theorem \ref{thm:alg-post-proc} is the direct consequence of Theorem \ref{thm:main-BIM}.
Algorithm \ref{Alg:post-proc} did not fail in any of our 
simulations in Section \ref{sec:simulation} of IEEE test cases. 

Similarly for BFM, one could also apply post-processing to recover the solution of \eqref{eq:OPF-BFM} from an optimal solution of \eqref{eq:ROPF-BFM}.
In the BFM, instead of checking the rank of $\W^*$, we check the rank of $\Mat^{\VolSqr^*,\PF,^*\CurSqr^*}(j,k)$ for each $j\sim k$ and $\VolSqr_j^*$ for $j\in\Vertex$.

\subsection{Relaxation with Penalized Cost Function}\label{subsec:penalized}
\revise{
Since the inexactness of relaxations \eqref{eq:ROPF-BIM} and \eqref{eq:ROPF-BFM} originates from two issues: the non-uniqueness in $\CurDeltaSqr^*$ and the spectrum error, where the latter is essentially amplified by the former.
The second algorithm we propose is to penalize and suppress the trace of $\CurDeltaSqr_j$ in the cost function.
With such penalty term, the value of $\CurDeltaSqr^*$ will be unique for fixed $\W^*$ and $\X^*$ in the solution of \eqref{eq:ROPF-BIM} and the spectrum error can also be restricted.} Similar penalization approaches were previously proposed in~\cite{madani2015promises,molzahn2016laplacian} to promote low-rank solutions of relaxed problems.

The penalized relaxed problem formulation under the BIM becomes 
\begin{subequations}
\begin{align}
\underset{\s,\s_\Delta,\W,\X,\CurDeltaSqr}{\text{minimize}} &&& f(\s,\s_{\Delta})+\lambda\sum\limits_{j\in\Vertex}\tr(\rho_j)
\label{eq:OPF-BIM-penalty.a}\\
\text{\quad subject to}&&&  \eqref{eq:PF-BIM.a}-\eqref{eq:PF-BIM.d}, \eqref{eq:PF-BIM.f},
\eqref{eq:OPF.d}, \eqref{eq:OPF-BIM.c}
\label{eq:OPF-BIM-penalty.b}
\end{align}
\label{eq:OPF-BIM-penalty}
\end{subequations}
Similarly, the penalized relaxed program under the BFM becomes
\begin{subequations}
\begin{align}
\underset{\s,\s_\Delta,\PF,\VolSqr,\CurSqr,\X,\CurDeltaSqr}{\text{minimize}} &&& f(\s,\s_{\Delta})+\lambda\sum\limits_{j\in\Vertex}\tr(\rho_j)
\label{eq:OPF-BFM-penalty.a}\\
\text{\quad subject to}&&& \eqref{eq:OPF.d}, \eqref{eq:OPF-BFM.c}, \eqref{eq:PF-BFM.a}-\eqref{eq:PF-BFM.g}, \eqref{eq:PF-BFM.i}.
\label{eq:OPF-BFM-penalty.b}
\end{align}
\label{eq:OPF-BFM-penalty}
\end{subequations}
Because $\tr(\CurDeltaSqr_j)$ is linear 
and all constraints in \eqref{eq:OPF-BIM-penalty.b} and \eqref{eq:OPF-BFM-penalty.b} are convex,
both \eqref{eq:OPF-BIM-penalty} and \eqref{eq:OPF-BFM-penalty} are convex optimization problems and can be efficiently solved in polynomial time.\footnote{If the constraint \eqref{eq:OPF.d} has a linear or quadratic formulation, then both problems are SOCPs.}
Here, $\lambda>0$ controls the weight of $\sum\tr(\rho_j)$ in the cost function.
The pseudo code (based on BIM) is summarized in Algorithm \ref{Alg:penalized}. The algorithm for BFM is similar.

\begin{algorithm}
\caption{Relaxation Algorithm with Penalized Cost Function.}
\begin{algorithmic}[1]\label{Alg:penalized}
 \renewcommand{\algorithmicrequire}{\textbf{Input:}}
 \renewcommand{\algorithmicensure}{\textbf{Output:}}
 \REQUIRE $y$, $\set$, $\set_{\Delta}$
 \ENSURE  Optimal solution $(\s,\s_\Delta,\W,\X,\CurDeltaSqr)$ to \eqref{eq:OPF-BIM}.
  \STATE Pick a sufficiently small $\lambda>0$
  \STATE Solve \eqref{eq:OPF-BIM-penalty} and obtain $\sol^*:=(\s^*,\s_{\Delta}^*,\W^*,\X^*,\CurDeltaSqr^*)$
  \IF {($\rank({\W}^*)>1$)}
  \STATE Output `Failed!'
  \STATE Exit
  \ELSE
  \RETURN $\sol^*$
  \ENDIF
\end{algorithmic} 
\end{algorithm} 

Because the cost function in the penalized program has been changed, the output of Algorithm \ref{Alg:penalized} might not be the global optimal solution of \eqref{eq:OPF-BIM}.
We next show that the output of Algorithm \ref{Alg:penalized} serves as an approximation of the true optimal solution.
We make the following assumption.
\begin{assumption}
The problem \eqref{eq:ROPF-BIM} has at least one finite optimal solution.
\end{assumption}

Now consider a sequence of positive and decreasing $\lambda_i$ for $i=1,2,\cdots$.
Taking BIM as an example,
let the optimal solution of \eqref{eq:OPF-BIM-penalty} with respect to $\lambda_i$ be $\sol^{(i)}$.\footnote{If the program has multiple solutions, then pick any one of them.}
Then the following lemma implies the sequence $\sol^{(i)}$ has a limit point.
\begin{lemma}
The sequence $(\sol^{(i)})_{i=1}^\infty$ resides in a compact set, and hence has a limit point.
\end{lemma}
\begin{IEEEproof}
Because all the constraints in \eqref{eq:OPF-BIM-penalty} are closed, we only need to prove boundedness. 
By assumption, $\s_j^{(i)}$ and $\s_{\Delta,j}^{(i)}$ at bus $j$ are in compact sets $\set_j$ and $\set_{\Delta,j}$ respectively. The positive semi-definite matrix $\W^{(i)}$ has upper bounds on its diagonal elements and is therefore bounded. We only need to show that $\sum_j\tr({\CurDeltaSqr}_j^{(i)})$ is also bounded because the boundedness of $\X^{(i)}$ is implied by the constraint \eqref{eq:PF-BIM.f} as long as $\sum_j\tr({\CurDeltaSqr}_j^{(i)})$ is bounded.

To show $\sum_j\tr({\CurDeltaSqr}_j^{(i)})$ is bounded, let $\hat{\sol}=(\hat{\s},\hat{\s}_{\Delta},\hat{\W},\hat{\X},\hat{\CurDeltaSqr})$ be an optimal solution of \eqref{eq:ROPF-BIM}.  Then $\hat{\sol}$ is feasible for \eqref{eq:OPF-BIM-penalty} regardless of the value of $\lambda$.
For any $i$, we must have $\sum_j\tr({\CurDeltaSqr}_j^{(i)})\leq\sum_j\tr(\hat{\CurDeltaSqr}_j)$; otherwise, $\hat{\sol}$ will always give a strictly smaller cost value in \eqref{eq:OPF-BIM-penalty} for $\lambda=\lambda_i$ and it would contradict the optimality of $\sol^{(i)}$.
\end{IEEEproof}

Suppose $\tilde{\sol}:=(\tilde{\s},\tilde{\s}_\Delta,\tilde{\W},\tilde{\X},\tilde{\CurDeltaSqr})$ is an arbitrary limit point of the sequence $\sol^{(i)}$. We present sufficient conditions for $\tilde{\sol}$ to be an optimal solution of~\eqref{eq:OPF-BIM}. First, we introduce the following lemma, which will be used to prove the optimality conditions.

\begin{lemma}\label{lem:mintrace}
Consider the positive semi-definite matrix $\Mat$ as in~\eqref{eq:block-psd} where $\AMat=\xVec\xVec^\Hn$ for some vector $\xVec$ such that $\xVec\neq\zero$, and $\BMat = \xVec\yVec^\Hn$. Then, $$\tr(\Mat)\geq\xVec^{\Hn}\xVec+\yVec^{\Hn}\yVec$$ and equality holds if and only if $\CMat=\yVec\yVec^\Hn$.
\end{lemma}
\begin{IEEEproof}
It is sufficient to prove that $\CMat-\yVec\yVec^\Hn\succeq 0$.
If not, then suppose there exists $\zVec$ such that $\zVec^\Hn(\CMat-\yVec\yVec^\Hn)\zVec<0$.
Because $\xVec\neq\zero$, we can always find $\wVec$ such that $\wVec^\Hn\xVec=-\zVec^\Hn\yVec$.
Consider
\begin{align*}
&\left[
    \begin{array}{c}
        \wVec  \\
        \zVec 
    \end{array}
\right]^\Hn\Mat
\left[
    \begin{array}{c}
        \wVec  \\
        \zVec 
    \end{array}
\right]=
\left[
    \begin{array}{c}
        \wVec  \\
        \zVec 
    \end{array}
\right]^\Hn
\left[
\begin{array}{cc}
    \xVec\xVec^\Hn & \xVec\yVec^\Hn \\
    \yVec\xVec^\Hn & \CMat
\end{array}
\right]
\left[
    \begin{array}{c}
        \wVec  \\
        \zVec 
    \end{array}
\right]\\
=& \ \wVec^\Hn\xVec\xVec^\Hn\wVec + \zVec^\Hn\yVec\xVec^\Hn\wVec
+ \wVec^\Hn\xVec\yVec^\Hn\zVec + \zVec^\Hn\CMat\zVec\\
=& \ \zVec^\Hn\yVec\yVec^\Hn\zVec-\zVec^\Hn\yVec\yVec^\Hn\zVec-\zVec^\Hn\yVec\yVec^\Hn\zVec+\zVec^\Hn\CMat\zVec\\
=& \ \zVec^\Hn(\CMat-\yVec\yVec^\Hn)\zVec<0.
\end{align*}
This contradicts the positive semi-definiteness of $\Mat$.
\end{IEEEproof}

\begin{theorem}\label{thm:seq}
If $\rank(\tilde{\W})=1$, then $\tilde{\sol}$ is an optimal solution of \eqref{eq:OPF-BIM}.
\end{theorem}

\begin{IEEEproof}
We first show that $\tilde{\sol}$ is the optimal solution of \eqref{eq:ROPF-BIM}. Since \eqref{eq:ROPF-BIM} and \eqref{eq:OPF-BIM-penalty} have the same feasible set, which is closed, $\tilde{\sol}$ is also feasible for \eqref{eq:ROPF-BIM} and \eqref{eq:OPF-BIM-penalty} for any $\lambda$. If $\tilde{\sol}$ is not optimal for \eqref{eq:ROPF-BIM}, then there must exist another point $\bar{\sol}$ such that $f(\bar{\s},\bar{\s}_{\Delta})+\alpha=f(\tilde{\s},\tilde{\s}_{\Delta})$ and $\alpha>0$. 
Then for some sufficiently large $i_0$, we have
\begin{align*}
\lambda_{i_0}\sum_{j\in\Vertex}\tr(\bar{\CurDeltaSqr}_j)<\frac{\alpha}{2}\\
|f(\tilde{\s},\tilde{\s}_\Delta)-f(\s^{(i_0)},\s_\Delta^{(i_0)})|<\frac{\alpha}{2}.
\end{align*}
Therefore:
\begin{align*}
f(\bar{\s},\bar{\s}_{\Delta})+\lambda_{i_0}\sum_{j\in\Vertex}\tr(\bar{\CurDeltaSqr}_j)<
f(\s^{(i_0)},\s_\Delta^{(i_0)})+\lambda_{i_0}\sum_{j\in\Vertex}\tr({\CurDeltaSqr}_j^{(i_0)})
\end{align*}
which contradicts the optimality of $\sol^{(i_0)}$.

Then for each $j$, we decompose $\tilde{\W}_{jj}=\tilde{\xVec}_j\tilde{\xVec}_j^\Hn$, $\tilde{\X}=\tilde{\xVec}_j\tilde{\yVec}_j^\Hn$, and construct ${\CurDeltaSqr}_j^\dagger$ as $\tilde{\yVec}_j\tilde{\yVec}_j^\Hn$.
Under the same argument as in the proof of Theorem \ref{thm:main-BIM},
the solution ${\sol}^\dagger:=(\tilde{\s},\tilde{\s}_{\Delta},\tilde{\W},\tilde{\X},{\CurDeltaSqr}^\dagger)$ is an optimal solution for both \eqref{eq:OPF-BIM} and \eqref{eq:ROPF-BIM}. 
We want to prove $\tilde{\sol}={\sol}^\dagger$ and conclude that $\tilde{\sol}$ is also an optimal solution for \eqref{eq:OPF-BIM}. The proof is by contradiction.

Since ${\sol}^\dagger$ is optimal for \eqref{eq:OPF-BIM}, $\Mat^{\tilde{\W},\tilde{\X},\CurDeltaSqr^\dagger}(j)$ must be of rank-1 for all $j$.
By Lemma \ref{lem:mintrace}, we have
\begin{align*}
\tr(\Mat^{\tilde{\W},\tilde{\X},\tilde{\CurDeltaSqr}}(j))
\geq
\tr(\Mat^{\tilde{\W},\tilde{\X},\CurDeltaSqr^\dagger}(j))
\end{align*}
and therefore $\tr(\tilde{\CurDeltaSqr}_j)\geq\tr({\CurDeltaSqr}_j^\dagger)$ for all $j$.
If $\tilde{\sol}\neq{\sol}^\dagger$, then some equalities cannot be achieved, and as a result, $\sum_j\tr(\tilde{\CurDeltaSqr}_j)-\sum_j\tr({\CurDeltaSqr}_j^\dagger)=\beta$ for some $\beta>0$.

As $\tilde{\sol}$ is a limit point of $(\sol^{(i)})_{i=1}^\infty$, there must be some sufficiently large $i_1$ such that
\begin{align*}
    \Big|\sum_{j\in\Vertex}\tr(\tilde{\CurDeltaSqr}_j)-\sum_{j\in\Vertex}\tr({\CurDeltaSqr}_j^{(i_1)})\Big|<\frac{\beta}{2},
\end{align*}
Hence
\begin{align*}
    \sum_{j\in\Vertex}\tr({\CurDeltaSqr}_j^\dagger)<\sum_{j\in\Vertex}\tr({\CurDeltaSqr}_j^{(i_1)}).
\end{align*}
On the other hand, \eqref{eq:ROPF-BIM} and \eqref{eq:OPF-BIM-penalty} have the same feasible set, so the optimality of ${\sol}^\dagger$ for \eqref{eq:ROPF-BIM} implies 
$f(\tilde{\s},\tilde{\s}_{\Delta})\leq f\left(\s^{(i_1)},{\s}_{\Delta}^{(i_1)} \right)$.  Therefore
\begin{align*}
    f(\tilde{\s},\tilde{\s}_{\Delta})+\lambda_{i_1}\sum_{j\in\Vertex}\tr({\CurDeltaSqr}_j^\dagger) < f\left(\s^{(i_1)},{\s}_{\Delta}^{(i_1)} \right)+\lambda_{i_1}\sum_{j\in\Vertex}\tr\left({\CurDeltaSqr}_j^{(i_1)}\right)
\end{align*}
which contradicts the fact that $\sol^{(i_1)}$ is the optimal solution for \eqref{eq:OPF-BIM-penalty} with respect to $\lambda_{i_1}$.
\end{IEEEproof}

Theorem \ref{thm:seq} shows that when we solve the penalized program with a sequence of decreasing $\lambda_i$, any limit point would be a global optimal for \eqref{eq:OPF-BIM} as long as the $\W$ matrix associated with the limit point is of rank-1. In our simulations, we apply Algorithm \ref{Alg:penalized} to solve \eqref{eq:OPF-BIM-penalty} with a fixed but sufficiently small $\lambda$, which usually results in rank-1 solutions.

\begin{remark}
Further, if all optimal solutions of \eqref{eq:ROPF-BIM} have the same value for $\s,\s_\Delta,\W,\X$, then Algorithm \ref{Alg:post-proc} succeeds if and only if $\rank(\tilde{\W})=1$ holds in Theorem \ref{thm:seq}. If Algorithm \ref{Alg:post-proc} succeeds, its output will also be the same as $\tilde{\sol}$.
\end{remark}
\section{Model Equivalence} \label{sec:analysis}
In the previous sections, our results for the BIM and BFM always come in pairs and are analogous. 
A natural question is whether there exist instances where
one model produces an exact solution while the other does not.
In single-phase networks and multi-phase systems with only wye connections, \cite{bose2014equivalent} and \cite{gan2014convex} have shown that the two models are equivalent in the sense that one will produce an exact solution if and only if the other will. 
We show in this subsection that a similar result holds in the presence of delta connections. 
We first define the equivalence between two optimization problems as follows.
\begin{definition}\label{def:equiv}
Consider two optimization problems
\begin{subequations}
\begin{align}
\underset{x}{\text{minimize}} &&& f_A(x)
\label{eq:optA.a}\\
\text{\quad subject to}&&& x\in\feasiblesetX
\label{eq:optA.b}
\end{align}
\label{eq:optA}
\end{subequations}
and
\begin{subequations}
\begin{align}
\underset{y}{\text{minimize}} &&& f_B(y)
\label{eq:optB.a}\\
\text{\quad subject to}&&& y\in\feasiblesetY.
\label{eq:optB.b}
\end{align}
\label{eq:optB}
\end{subequations}
We say \eqref{eq:optA} and \eqref{eq:optB} are {\it equivalent} if there exist mappings $g_1:\feasiblesetX\rightarrow\feasiblesetY$ and $g_2:\feasiblesetY\rightarrow\feasiblesetX$
such that 
\begin{align*}
x\in\feasiblesetX \Rightarrow g_1(x)\in\feasiblesetY,~f_A(x)=f_B(g_1(x)),\\
y\in\feasiblesetY \Rightarrow g_2(y)\in\feasiblesetX,~f_B(y)=f_A(g_2(y)).
\end{align*}
\end{definition}
We do not require $g_1$ and $g_2$ to be bijections, but if one of the mappings is a bijection, then we can always set the other as its inverse.
We denote the decision variables for the BIM as
\begin{align*}
\sol^\BIM=(\s^\BIM,\s_\Delta^\BIM,\W^\BIM,\X^\BIM,\CurDeltaSqr^\BIM)
\end{align*}
and the decision variables for the BFM as
\begin{align*}
\sol^\BFM=(\s^\BFM,\s_\Delta^\BFM,\PF^\BFM,\VolSqr^\BFM,\CurSqr^\BFM,\X^\BFM,\CurDeltaSqr^\BFM).
\end{align*}
The superscripts here are to distinguish the same variable for different models. 

\begin{proposition}\label{prop:equiv-R}
Problems \eqref{eq:ROPF-BIM} and \eqref{eq:ROPF-BFM} are equivalent.
Moreover, for the pairs $g_1$ and $g_2$ in Definition \ref{def:equiv}, if $\sol^{\BIM}$ satisfies \eqref{eq:PF-BIM.e}, then $g_1(\sol^\BIM)$  satisfies \eqref{eq:PF-BFM.h}. If $\sol^{\BFM}$ satisfies \eqref{eq:PF-BFM.h}, then $g_2(\sol^\BFM)$  satisfies \eqref{eq:PF-BIM.e}. 
\end{proposition}

We only sketch a proof here by providing the mappings $g_1$ and $g_2$, where $g_1$ can be written as
\begin{subequations}
\begin{align}
\s^\BFM &= \s^\BIM,~\s_\Delta^\BFM = \s_\Delta^\BIM\\
\PF^\BFM_{jk} &= (\W_{jj}^\BIM-\W_{jk}^\BIM) y_{jk}^\Hn\\
\VolSqr_j^\BFM &= \W_{jj}^\BIM\\
\CurSqr_{jk}^\BFM &= y_{jk}(\W_{jj}^\BIM\!+\!\W_{kk}^\BIM\!-\!\W_{jk}^\BIM\!-\!\W_{kj}^\BIM) y_{jk}^\Hn\\
\X_j^\BFM &= \X_j^\BIM,~\CurDeltaSqr_j^\BFM = \CurDeltaSqr_j^\BIM
\end{align}\label{eq:equiv-BIM-BFM}
\end{subequations}
and $g_2$ as
\begin{subequations}
\begin{align}
\s^\BIM &= \s^\BFM,~\s_\Delta^\BIM = \s_\Delta^\BFM\\
\W_{jj}^\BIM &= \VolSqr_j^\BFM\\
\W_{jk}^\BIM &=
\left\{
\begin{array}{ll}
    \VolSqr_j^\BFM-\PF_{jk}^\BFM z_{jk}^\Hn, & \text{if}~j\rightarrow k  \\
    (\W_{kj}^\BIM)^\Hn, & \text{if}~k\rightarrow j  \\
\end{array}
\right.
\\
\X_j^\BIM &= \X_j^\BFM,~\CurDeltaSqr_j^\BIM = \CurDeltaSqr_j^\BFM.
\end{align}\label{eq:equiv-BFM-BIM}
\end{subequations}
For $g_2$, the value of $\W_{jk}^\BIM$ where $j\neq k$ and $(j,k)\not\in\Edge$ can be determined arbitrarily as long as $\W\succeq 0$. As $\Graph$ is a tree, the values in \eqref{eq:equiv-BFM-BIM} can guarantee that we can always complete the matrix $\W^\BIM\succeq 0$, but not necessarily in a unique way.

Proposition \ref{prop:equiv-R} shows that to apply Algorithm \ref{Alg:post-proc}, if an optimal solution of \eqref{eq:ROPF-BIM} can produce an exact solution of \eqref{eq:OPF-BIM}, then there must also be an optimal solution of \eqref{eq:ROPF-BFM} that can produce an exact solution of \eqref{eq:OPF-BFM}, 
\update{even though both \eqref{eq:OPF-BIM} and \eqref{eq:OPF-BFM} may have multiple solutions.}
The converse is also true. Informally, for Algorithm \ref{Alg:post-proc}, both the BIM and BFM have the same capability of producing exact solutions.

The same holds for the penalized program.

\begin{proposition}
Problems \eqref{eq:OPF-BIM-penalty} in the BIM and \eqref{eq:OPF-BFM-penalty} in the BFM are equivalent when $\lambda$ takes the same value for both problems.
\end{proposition}

The proposition can be easily proved using the same mappings $g_1$ and $g_2$ in \eqref{eq:equiv-BIM-BFM} and \eqref{eq:equiv-BFM-BIM}, respectively.

\update{Beyond OPF problems, mappings $g_1$ and $g_2$ also provide the correspondence between feasible points under the two models.
Thereby, a solution of power flow equations under one model can also be translated into a solution with the same physical meaning under the other model by applying $g_1$ or $g_2$. Note that power flow equations may have multiple solutions.}


\update{Though BIM and BFM are mathematically equivalent, the two models may behave differently in practice and shed lights on different properties.
Some analysis may rely on the structure of one model but not the other, which is indeed the case for single-phase networks \cite{low2014convexI}. The equivalence implies that one could freely choose a model that is more convenient for a specific problem. For instance, a recent work \cite{zhou2019sufficient} on exact relaxation for three-phase networks was derived based on BIM.}
\section{Numerical Results} \label{sec:simulation}
In this section, we show the ability of the proposed relaxation algorithms to recover the optimal solution to~\eqref{eq:OPF-BIM} and~\eqref{eq:OPF-BFM}. We use the IEEE 13-, 37-, \update{and 123}-node distribution feeders \cite{schneider2017analytic} to assess the exactness of both algorithms for both the BIM and BFM models. \update{Note that the IEEE 123-bus feeder does not include delta-connected components. Hence, we artificially added $4$ delta-connected loads to the feeder to assess the performance of the proposed approaches.}
In our experiments, we check how close the output matrices $\W, \MWXrho, \MvSl, \MvXrho$ are to being rank-1, and we evaluate the maximum violation of the constraints when the decision variables are produced from the two proposed algorithms. For all the experiments in this section, we show that both algorithms succeed up to numerical precision, and each has its own advantages and disadvantages.

\subsection{Experimental Setup}
The load transformer in the IEEE test feeders are modeled as lines with equivalent impedance, whereas the substation transformers and regulators are removed. The switches are assumed to be open or short according to their default status. The capacitor banks are modeled as controllable reactive power sources with continuous control space.

The voltage at the substation is assumed to be $ \Vref = \overline{V}[1, e^{\sfrac{-\iu2\pi}{3}}, e^{\sfrac{\iu2\pi}{3}}]^\T$, where $\overline{V}$ is the maximum allowed voltage magnitude.
\update{
The operational constraints for controllable loads are set as in \cite{zhao2017optimal}.
}
The AC-OPF problem is solved with the cost function $f(\s, \s_\Delta)$ comprising three parts. The first part minimizes the total power losses in the network, and it can be written as
\begin{align*}
    p_{\text{loss}} = \sum_{j \in \Vertex} \sum_{\phi \in \phases_\Y^{j}} \pya{j}^{\phi} + \sum_{j \in \Vertex} \sum_{\phi \in \phases_\Delta^{j}} \pda{j}^{\phi}.
\end{align*}

The second part penalizes deviations of  the active and reactive injection profile from nominal profiles, and it is given by
$d_p(\py, \pd)$ and $d_q(\qy, \qd)$ as follows.
\begin{align}\notag
    d_p(\py, \pd) = &\sum_{j \in \Vertex} \sum_{\phi \in \phases_\Y^{j}} \frac{1}{2\pyref} ( \pya{j}^{\phi} - \pyref)^2 \\ \notag
    &+ \sum_{j \in \Vertex} \sum_{\phi \in \phases_\Delta^{j}} \frac{1}{2\pdref} ( \pda{j}^{\phi} - \pdref)^2,
\end{align}
\begin{align}\notag
    d_q(\qy, \qd) = &\sum_{j \in \Vertex} \sum_{\phi \in \phases_\Y^{j}} \frac{1}{2\qyref} ( \qya{j}^{\phi} - \qyref)^2 \\ \notag
    &+ \sum_{j \in \Vertex} \sum_{\phi \in \phases_\Delta^{j}} \frac{1}{2\qdref} ( \qda{j}^{\phi} - \qdref)^2.
\end{align}
The values $\pyref,\qyref,\pdref,\qdref$ represent the nominal active and reactive injection values for phase $\phi$ at bus $j$. All the tracking errors are normalized by their nominal values to have the same order of magnitude for all quantities.
In addition, ${\phases}_{\Y}^j \subseteq \{a, b, c\}$ and ${\phases}_\Delta^j \subseteq \{ab, bc, ca\}$ denote the available wye and delta connections at bus $j \in \Vertex$, respectively. 
\update{Penalizing the deviation of power injection can characterize either the operational cost of controllable loads, the curtailment of photovoltaic systems, or the charging/discharging cost of batteries. The same cost expression was also used in \cite{zhao2017optimal, dall2017optimal}.}

The last part minimizes the deviation of the power injections at the substation from the reference injections $\overline{p}_0, \overline{q}_0\in\Real$ provided by the transmission system operator. Therefore, the system operational cost function can be written as
\begin{align}\notag
    f(\s, \s_\Delta) = &  \mu_\ell\ p_{\text{loss}} + w_{p}\ d_{p}(\py, \pd) + w_{q}\ d_{q}(\qy, \qd) \\ & +  \mu_p\ \frac{(\one^\T \py_{0} - \overline{p}_{0})^2}{\overline{p}_0} + \mu_q\ \frac{(\one^\T \qy_{0} - \overline{q}_{0})^2}{\overline{q}_0}. \notag
\end{align}
The nonnegative weights $w_p$, $w_q$, $ \mu_\ell$, $\mu_p$ and $\mu_q$ are used to reflect the relative importance of the components of the cost function and are set as follows
\begin{align}
    w_p = w_q = \mu_\ell = 1, \qquad \mu_p = \mu_q = 4.
\end{align}




\begin{table*}
	\renewcommand{\arraystretch}{1.1}
	\caption[]{Rank and infeasibility for the outputs of Algorithm \ref{Alg:post-proc} (with post-processing).}
	\begin{center}
		\begin{tabular}{l  c  c c  c c} 
			\toprule
			\multirow{2}{*}{Network} &  \multirow{2}{*}{Voltage} & \multicolumn{2}{c}{BIM} & \multicolumn{2}{c}{BFM}\\
			\cmidrule(lr){3-4}
            \cmidrule(lr){5-6}
			&&{$\W$-ratio}  & {Infeas. (kW)}&{$\MvSl$-ratio} &  {Infeas. (kW)}\\
			\midrule
			\multirow{2}{*}{\em IEEE-13} &  $3\%$ & $ 1.91 \times 10^{-8}$ & $3.20 \times 10^{-1} $ & $1.74\times 10^{-5} $& $4.29 \times 10^{-2}$  \\ 
			& $5\%$ & $2.73 \times 10^{-9}$ & $3.20 \times 10^{-1}$ & $1.58\times 10^{-5} $& $4.10\times 10^{-2}$\\
			\midrule
			\multirow{2}{*}{\em IEEE-37} &  $3\%$ & $ 4.81 \times 10^{-10}$ & $9.84 \times 10^{-2}$ & $7.81 \times 10^{-5} $ & $9.67\times10^{-2}$ \\ 
			& $5\%$ & $2.77 \times 10^{-8}$ & $9.83\times10^{-2}$ & $2.70\times 10^{-5}$  & $9.72\times10^{-2}$\\
			\midrule
			\multirow{2}{*}{\em IEEE-123} &  $3\%$ & $ 7.75 \times 10^{-8}$ & $1.54 \times 10^{-3}$ & $1.07 \times 10^{-4} $  & $1.03\times10^{-2}$ \\ 
			& $5\%$ & $7.67 \times 10^{-8}$ & $1.54\times10^{-3}$ & $9.64\times 10^{-5}$ &  $1.02\times10^{-2}$\\
			\bottomrule
		\end{tabular}
	\end{center}
	\label{table:alg-post-proc}
\end{table*}

\begin{table*}
	\renewcommand{\arraystretch}{1.1}
	\caption[]{Rank and infeasibility for the outputs of Algorithm \ref{Alg:penalized} (with penalized cost function).}\vspace{-10pt}
	\begin{center}
		\begin{tabular}{l c c c c c c c} 
			\toprule
			\multirow{2}{*}{Network} &  \multirow{3}{*}{Voltage} &   \multicolumn{3}{c}{BIM} & \multicolumn{2}{c}{BFM} \\
			\cmidrule(lr){3-5}
            \cmidrule(lr){6-8}
			& & {$ \W$-ratio} & {$\MWXrho$-ratio} & {Infeas. (kW)} & {$\MvSl$}-ratio &{$\MvXrho$-ratio} & {Infeas. (kW)}\\
			\midrule
			\multirow{2}{*}{\em IEEE-13} &  $3\%$ & $1.34 \times 10^{-10}$  & $2.36\times 10^{-9}$ & $8.85 \times 10^{-2}$ & $1.44 \times 10^{-10}$ & $1.97 \times 10^{-10}$ & $4.43 \times 10^{-5}$\\
			& $5\%$ & $1.31 \times 10^{-10}$ & $1.96 \times 10^{-9}$ & $8.84 \times 10^{-2}$ &$1.36 \times 10^{-10}$ & $1.57 \times 10^{-10}$ & $1.46 \times 10^{-5}$\\
			\midrule
			\multirow{2}{*}{\em IEEE-37} &  $3\%$  & $3.04\times 10^{-8} $&$6.22 \times 10^{-8} $ &$5.75 \times 10^{-6}$ & $8.85 \times 10^{-8} $ &$3.38 \times 10^{-5} $ &$1.45 \times 10^{-6}$ \\ 
			& $5\%$ & $ 2.94 \times 10^{-8}$ &$1.05\times 10^{-8}$ & $1.06\times 10^{-6}$ & $2.12 \times 10^{-8} $&$3.18\times 10^{-5}$ & $1.00\times 10^{-6}$\\
			\midrule
			\multirow{2}{*}{\em IEEE-123} &  $3\%$  & $1.45\times 10^{-9} $&$7.03 \times 10^{-9} $ &$9.13 \times 10^{-7}$ & $1.05 \times 10^{-8} $ &$8.99 \times 10^{-9} $ &$1.34 \times 10^{-6}$ \\ 
			& $5\%$ & $ 1.94 \times 10^{-8}$ &$9.31\times 10^{-8}$ & $7.84\times 10^{-6}$ & $7.98 \times 10^{-9} $&$6.59\times 10^{-9}$ & $1.40\times 10^{-6}$\\
			\bottomrule
		\end{tabular}
	\end{center}
	\label{table:alg-penalized}\vspace{-10pt}
\end{table*}

\begin{table}
	\renewcommand{\arraystretch}{1.1}
	\caption[]{Effect of the penalty parameter on the cost and infeasibility.}\vspace{-20pt}
	\begin{center}
		\begin{tabular}{l c c c c} 
			\toprule
			\multirow{2}{*}{$\boldsymbol{\lambda}$} &  \multicolumn{2}{c}{BIM} &   \multicolumn{2}{c}{BFM} \\
			\cmidrule(lr){2-3}
            \cmidrule(lr){4-5}
			&  Cost  & {Infeas. (kW)} & Cost  & {Infeas. (kW)}\\
			\midrule
			$0$ &  $100.0036$ & $9.84 \times 10^{-2}$ & $100.0194$ & $ 9.67 \times 10^{-2}$\\
			$0.1$& $103.9504$ & $1.15 \times 10^{-2}$ & $104.7141$ & $ 5.99 \times 10^{-4}$\\
			$1$ &  $104.7846$ & $6.00 \times 10^{-5}$ & $104.7840$ & $ 1.80 \times 10^{-5}$\\
			$10$&  $104.7886$ & $5.75 \times 10^{-6}$ & $104.7982$ & $ 1.45 \times 10^{-6}$\\
			$100$& $104.9431$ & $3.17 \times 10^{-6}$ & $105.0332$ & $ 1.23 \times 10^{-7}$\\
			\bottomrule
		\end{tabular}
	\end{center}
	\label{table:penalty-effect}\vspace{-10pt}
\end{table}

\subsection{Exactness Results for Algorithm~\ref{Alg:post-proc}}
In this subsection, we assess the quality of the solutions recovered using Algorithm \ref{Alg:post-proc}. We solve \eqref{eq:ROPF-BIM} for the BIM as well as \eqref{eq:ROPF-BFM} for the BFM with different values of voltage limits for the three considered feeders. We invoke the Mosek conic solver using CVX, a MATLAB-based convex optimization toolbox. 

The left-hand side of Table~\ref{table:alg-post-proc} provides the result of Algorithm \ref{Alg:post-proc} based on the BIM. The voltage column represents the maximum and minimum voltage deviation allowed, i.e., $3\%$ means that the value of $ \overline{V} $ and $\underline{V}$ are set to $1.03$ and ${0.97}$, respectively. We assess the rank of matrices $\W_{jj}$, for all $j\in \Vertex$, in terms of the ratio between the top two largest eigenvalues of these matrices. The maximum ratio among all $j\in\Vertex$ is listed in the table. \revise{In the solution of~\eqref{eq:ROPF-BIM} (before post-processing), the ratio between the two maximum eigenvalues of the matrices $\MWXrho$ is on the order of $10^{-1}$, and after the post-processing in Algorithm \ref{Alg:post-proc}, the final $\MWXrho$-ratio will be dominated by $\W$-ratio and is thus not informative to be displayed in the table.}
Because of the the spectrum error, the output $\tilde{\X}$ could be different from $\X^*$, and thus having a very small $\W$-ratio is not enough to guarantee the feasibility of the final output of Algorithm~\ref{Alg:post-proc}.
Therefore, we also assess the infeasibilty of the power flow equations by measuring the maximum violation in \eqref{eq:PF-BIM.a} for the solutions returned by Algorithm~\ref{Alg:post-proc}.
Here, the violation is defined as the difference between the left- and right-hand sides of \eqref{eq:PF-BIM.a} (in kW) when $\s,\W,\X$ are evaluated as the output of Algorithm~\ref{Alg:post-proc}.

As Algorithm \ref{Alg:post-proc} circumvents the numerical error to a certain extent by recovering $\Cur_{\Delta}$ from~\eqref{eq:pf.b}, the infeasibility is on the order of $10^{-3}$ to $10^{-1}$ kW, which reflects the effect of the spectrum error.
\update{As a benchmark, the load injections for those feeders are on the order of $10^1$ to $10^2$ kW, and are typically two orders of magnitude higher than the infeasibility.}

On the right-hand side of Table~\ref{table:alg-post-proc}, the rank of $\MvSl(j,k)$ for all lines $(j,k)\in\Edge$ is examined for the same algorithm under the BFM. Again, we present the maximum ratio between the two largest eigenvalues. Similar to the BIM, the infeasibility, i.e., the violation of \eqref{eq:pf.a}, is shown in the table.



\subsection{Exactness Results for Algorithm \ref{Alg:penalized}}
In our setting, the penalized formulations~\eqref{eq:OPF-BIM-penalty} and~\eqref{eq:OPF-BFM-penalty} are solved with the parameter $\lambda = 10$ in all experiments. We will later show how the value of $\lambda$ affects the solution quality.

Table~\ref{table:alg-penalized} presents the maximum ratio between the top two largest eigenvalues of matrices $\MWXrho(j)$ for the BIM and $\MvXrho(j)$ for the BFM returned by the solvers. Comparing the infeasibility of the solutions obtained using Algorithm \ref{Alg:post-proc}, shown in Table \ref{table:alg-post-proc}, and Algorithm \ref{Alg:penalized}, shown in Table \ref{table:alg-penalized}, it is clear that adding a penalty helps reduce the effect of the spectrum error 
and leads to global optimal solutions with much lower infeasibility.



To assess the effect of the penalization approach on the quality of the solutions in terms of cost and feasibility, Table~\ref{table:penalty-effect} shows the effect of increasing the penalty parameter in the cost function as well as the maximum infeasibility of the power equations (in kVA) for the IEEE 37-bus network with $3\%$ voltage limits. 
\update{Note that the case $\lambda=0$ corresponds to the output of Algorithm \ref{Alg:post-proc}.}
Although the solution feasibility is enhanced by increasing the penalty parameter, the cost associated with the solution obtained also increases. Note that the cost obtained with no penalty, i.e., from Algorithm \ref{Alg:post-proc}, represents a lower bound for the optimal cost of the original AC-OPF problem. In addition, increasing the penalty parameter beyond the values considered in Table~\ref{table:penalty-effect} leads to uninteresting solutions because the cost function becomes dominated by the penalty term. 

\subsection{Algorithm Summary and Comparison}
\update{Algorithms \ref{Alg:post-proc} and \ref{Alg:penalized} can be useful for different applications. Algorithm \ref{Alg:post-proc} solves the un-penalized problem and therefore prioritizes cost minimization at the cost of larger constraint violation. This is because the recovered $\tilde{\X}$ could be different from $\X^*$ and therefore may not precisely satisfy the constraints such as \eqref{eq:PF-BIM.a}. That is the reason we evaluate the infeasibility as a main metric in Tables \ref{table:alg-post-proc} and \ref{table:alg-penalized}. The simulation shows that the infeasibility is typically two orders of magnitude smaller than the load injections in the network and should be acceptable. Algorithm \ref{Alg:penalized}, on the other hand, can recover an OPF solution with much smaller constraint violation, but the optimal cost is higher because of the penalty term.

We also benchmark the computational time of the proposed algorithms for solving the AC-OPF problems in our case studies. Since both algorithms require solving similar optimization problems with different cost functions, the computational time of both algorithms is similar. Hence, we only present the computational time of Algorithm 1 to solve the AC-OPF problem for all networks using both the BIM and BFM formulations. The algorithm was implemented using Mosek as a conic solver on a laptop with Intel Core i9 CPU (2.40 GHz), 16 GB RAM, macOS Catalina OS, and MATLAB R2019b. The results show that the proposed algorithms take less than $10$ seconds to solve the AC-OPF problem for the IEEE 123-bus network on a standard laptop, which demonstrates the computational efficiency of the proposed algorithms. Note that sparse semidefinite programming solvers, e.g.,~\cite{zhang2018sparse}, are expected to be very efficient in solving problems with thousands of buses.

\begin{table}
	\renewcommand{\arraystretch}{1.1}
	\caption[]{Computational Time for both BIM and BFM Formulations (in seconds).}
	\begin{center}
		\begin{tabular}{l c c c} 
			\toprule
			Model &  IEEE-13 & IEEE-37 & IEEE-123 \\
			\midrule
			BIM &  $2.46$ & $5.56$ & $9.03$ \\
			BFM &  $2.77$ & $4.93$ & $9.87$ 	\\		
			\bottomrule
		\end{tabular}
	\end{center}
	\label{table:comp-time}
\end{table}
}
\section{Discussion and Future Work}
\update{We now discuss the limitations and future directions of our work.
First, our result guarantees the exactness when $\W$ or $\MvSl$ is of rank 1. 
When $\W$ or $\MvSl$ is not of rank 1, our algorithms will fail to output feasible solutions but can still provide a lower bound for the optimal value.
Reference \cite{zhou2019sufficient} provides a sufficient condition to guarantee that $\W$ is of rank 1 for unbalanced systems without delta connections. It would be interesting to extend the result to networks with delta connections.
For single-phase networks, it is also known that when SDP or SOCP relaxation fails to provide a feasible solution, a heuristic post-processing can be applied to encourage feasibility. 
Another direction is whether the same techniques can be applied when the system is unbalanced and contains delta connections.

Second, the main focus of this paper is on relaxing the non-convexity due to quadratic power flow equations. There are  other sources of non-convexity in OPF problems. For instance, for systems with controllable switches and voltage regulators, there are additional decision variables which are discrete, and the underlying non-convexity is beyond the scope of this paper.  
The modelling of those devices has been studied in \cite{bazrafshan2019optimal,chiang1990optimal}.
Including switches and voltage regulators in the formulation of three-phase networks with delta connections is an important future work for more applications in distribution networks.

Finally, there is a trade-off between relaxation-based algorithms and local search algorithms for general non-convex programs. On one hand, local search algorithms are faster and more scalable for large systems compared to SDP relaxation. On the other hand, there is no guarantee that a local search algorithm will converge, and when it does the output is in general a local optimum.
There are also many active lines of work on improving the computational efficiency of convex relaxations, and combining the relaxation methods with local algorithms. More discussion could be found in the survey \cite{molzahn2019survey} and references therein.
}
\section{Conclusion}
\revise{This paper studied the SDP relaxation of the AC-OPF problem for an unbalanced three-phase radial network with delta connections, formulated under both the BIM and BFM. we showed the equivalence between the BIM and BFM formulations and presented sufficient conditions for recovering exact solutions of the nonconvex AC-OPF formulations from their respective relaxations. 
The paper also showed why conventional relaxation (by directly dropping rank-1 constraints) always fail even if the sufficient conditions we proposed are satisfied.
It is due to the non-uniqueness in the relaxation solution and the spectrum error in computation.
Inspired by this finding, we then proposed two algorithms which are guaranteed to produce exact solutions whenever our sufficient conditions are satisfied. One applies post-processing and produces lower cost but larger constraint violation. The other adds a penalty term and produces higher cost but smaller constraint violation.
In simulations, we demonstrated that for three IEEE standard test cases, both algorithms 
are able to recover near global optimal solutions with tolerable constraint violation and cost sub-optimality.}








\bibliographystyle{IEEEtran}
\bibliography{my-bibliography}
\end{document}